\documentclass[12pt]{article}

\newif\ifskiptext\skiptextfalse

\newif\ifomittext\omittexttrue

\makeatletter%shorter space after section nr.
\def\@seccntformat#1{\csname the#1\endcsname.\hspace*{0.5em}}

\def\@maketitle{%remove date
  \newpage
  \null
  \vskip 2em% 
  \begin{center}%
  \let \footnote \thanks
    {\LARGE \@title \par}%
    \vskip 1.5em%
    {\large
      \lineskip .5em%
      \begin{tabular}[t]{c}%
        \@author
      \end{tabular}\par}%
  \end{center}%
  \par
  \vskip 1.5em}%\@maketitle
\makeatother

\usepackage{natbib}
\def\citet{\citep} %%% no "text" citations in crossrefs in biblio

\def\citeP{\citep*} %%%ek3 * ==> full author list
\bibpunct[, ]{[}{]}{;}{a}{}{,}
\newcommand{\EKhref}[2]{URL: \url{#1}}
\usepackage[%
breaklinks%
,colorlinks%
,linkcolor=red% internal document links
,anchorcolor=blue%
,citecolor=blue%
,bookmarks=false%
]{hyperref}
\usepackage{breakurl}

\usepackage{amssymb}
\usepackage{amsmath}
\makeatletter
\newcommand{\tpmod}[1]{{\@displayfalse\pmod{#1}}}
\makeatother
\usepackage{mathdots}
\usepackage{longtable}

\usepackage{url} 
\usepackage{xspace}

\usepackage{theorem}

\global
\global

\makeatletter
\gdef\th@plain{\normalfont\itshape
  \def\@begintheorem##1##2{%
        \item[\hskip\labelsep \theorem@headerfont ##1\ ##2.]}%EK add .
\def\@opargbegintheorem##1##2##3{%
   \item[\hskip\labelsep \theorem@headerfont ##1\ ##2\ (##3).]}}
\makeatother

\newtheorem{theorem}{Theorem}[section]

{Algorithm}[section]%ek2 add

{\theorembodyfont{\upshape}

}%\theorembodyfont etc.

{\theorembodyfont{\upshape}
{Example}[section]%ek2 add
}%\theorembodyfont

{Problem}[section]%ek2 add

{\theoremstyle{plain}

}%theoremstyle etc.

\makeatletter
\def\@yproof[#1]{\@proof{ #1}}
\def\@proof#1{\begin{trivlist}\item[]{\em Proof#1.}}
\newenvironment{proof}{\@ifnextchar[{\@yproof}{\@proof{}
}}{~$\Box$\end{trivlist}}
\makeatother

\usepackage{enumerate}
\makeatletter
\def\@@enum@[#1]{%
  \@enLab{}\let\@enThe\@enQmark
  \@enloop#1\@enum@
  \ifx\@enThe\@enQmark\@warning{The counter will not be printed.%
   ^^J\space\@spaces\@spaces\@spaces The label is: \the\@enLab}\fi
  \expandafter\edef\csname label\@enumctr\endcsname{\the\@enLab}%
  \expandafter\let\csname the\@enumctr\endcsname\@enThe
  \csname c@\@enumctr\endcsname7
  \@enum@}
\makeatother

\newlength{\mylabelsep}

\newsavebox{\onedotspacemysep}\savebox{\onedotspacemysep}{1.{}\hskip\mylabelsep}
\newsavebox{\onelparpmysep}\savebox{\onelparpmysep}{1(a)\hskip\mylabelsep}
\newsavebox{\onelparpiidotmysep}\savebox{\onelparpiidotmysep}{1(a)ii.{}\hskip\mylabelsep}

\newcommand{\ruleaboveheader}[1]%lift in *ex
{\parbox[t][0cm][t]{0cm}{%
\hspace*{0cm}% horizontal positioning
\parbox[b][0cm][t]{0cm}{%
\vspace*{#1}% vertical positioning
\rule{\textwidth}{0.5pt}
\par}\par}\par\nobreak\vspace{-2ex}\noindent}

\newcommand{\rulebelowtext}[1]%depth
{\hfill
\parbox[t][0cm][t]{0cm}{%
\hspace*{-\textwidth}% horizontal positioning
\parbox[b][0cm][t]{0cm}{%
\vspace*{#1}% vertical positioning
\rule{\textwidth}{0.5pt}
\par}\par}}

\newcommand{\ourtriang}{\raise2.5pt\hbox{\tiny$\blacktriangleright$\ }}

\newcommand{\ZZ}{ {\mathbb Z}}
\newcommand{\QQ}{ {\mathbb Q}}

\newcommand{\KK}{ {\mathsf K}}
\newcommand{\LL}{ {\mathsf L}}
\newcommand{\NN}{ {\mathsf N}}
\newcommand{\FF}[1]{{\mathbb F}_{\!#1}}

\newcommand{\defequal}{%
\text{\vbox to 2.75ex{\vss\hbox{%
\begin{tabular}[b]{@{}c@{}}
\hbox to 0pt{\hss\footnotesize
\upshape def\hss}
\\[-1.5ex]
${}={}$
\end{tabular}}}}
}% newcommand

\title{A note on the
van der Waerden conjecture on
\\
random polynomials with symmetric Galois group
\\
for function fields
}%title

\author{%
Erich L. Kaltofen %${}^1$
\\
\\% empty line
\small %\llap{${}^1$}%
Dept.\ of Mathematics, North Carolina State University
\\[-0.5ex]
\small
Dept.\ of Computer Science, Duke University
\\[-0.5ex]
\small {\ttfamily kaltofen@ncsu.edu}; \url{https://users.cs.duke.edu/~elk27/}
}%author %>>>

\begin{document}

\maketitle

\begin{abstract}
\noindent
Let $f(x) = x^n + (a_{n-1}t + b_{n-1}) x^{n-1} + \cdots + (a_0 t + b_0)$ 
be of constant degree $n$ in $x$ and degree $\le 1$ in $t$, where
all $a_i,b_i$ are randomly and uniformly selected from a finite field
$\FF q$ of $q$ elements.  Then the probability that the Galois group
of $f$ over $\FF q(t)$ is the symmetric group $S_n$ on $n$ elements is
$1 - O(1/q)$.  Furthermore, the probability that the Galois group
of $f(x)$ over $\FF q(t)$ is not $S_n$ is $\ge 1/q$ for $n\ge 3$
and $> 1/q-1/(2q^2)$ for $n=2$.
\end{abstract}

%\begin{center}
%\date{\today}
%\end{center}

%\input{1intro}
\section{Introduction}

Let $n$ be an integer constant $\ge 1$.
The van der Waerden conjecture states that
\begin{multline}
\label{eq:vdW}
\text{Prob}
\Big( \text{Galois group over $\QQ$ of }
x^n + \sum_{i=0}^{n-1} a_i x^i
\text{ is not the symmetric group $S_n$}
\\
\big|\;a_i\in \ZZ\text{ and } | a_i | \le H
\Big) = O\big(\frac{1}{H}\big).
\end{multline}
In \citeP{bhargava2021galois} there is the first proof and 
an extensive bibliography;
see also \citeP{AGOL-DSZ21}{}.
Note that the probability (\ref{eq:vdW}) is asymptotically sharp:
for $a_0=0$ all polynomials are reducible and have a smaller Galois group.
From (\ref{eq:vdW}) one can derive (see Section~\ref{sec:rems})
the following function field analog:
\begin{multline}
\label{eq:vdWt}
\text{Prob}
\Big( \text{Galois group over $\QQ(t)$ of }
x^n + \sum_{i=0}^{n-1} (a_i t + b_i) x^i
\text{ is not the symmetric group $S_n$}
\\
\big|\;a_i,b_i\in \ZZ\text{ and } | a_i | \le H, | b_i | \le H
\Big) = O\big(\frac{1}{H^2}\big).
\end{multline}
Again, the probability (\ref{eq:vdWt}) is asymptotically sharp.
Here we consider the coefficient field $\FF q(t)$ where $\FF q$ is
a finite field with $q = p^\ell$ elements, for $\ell \ge 1$
and the prime characteristic $p \ge 2$.  We prove that
\begin{multline}
\label{eq:vdWq}
\text{Prob}
\Big( \text{Galois group over $\FF q(t)$ of }
x^n + \sum_{i=0}^{n-1} (a_i t + b_i) x^i
\text{ is not the symmetric group $S_n$}
\\
\big|\; a_i,b_i\in \FF q
\Big) = O\big(\frac{1}{q}\big).
\end{multline}
Again the probability (\ref{eq:vdWq}) is
asymptotically sharp: the probability that
$\text{GCD}(A(x), B(x)) \ne 1$ for $A(x) = \sum_{i=0}^{n-1} a_i x^i$
and $B(x) = x^n + \sum_{i=0}^{n-1} b_i x^i$ is exactly $1/q$
(see, for instance, \citeP{BB07}), so at least
$q^{2n-1}$ polynomials $x^n + \sum_{i=0}^{n-1} (a_i t + b_i) x^i$
have smaller Galois group.  Note that the Galois group of
a polynomial of degree~$n \ge 3$ in $\FF q[x]$ is not the symmetric group~$S_n$.
For $n=2$, one subtracts the $(q^2-q)/2$ irreducible
$x^2 + b_1 x + b_0$ at $a_1=a_0=0$ from the count.

%omittext

\section{Proof of Probability Estimate (\ref{eq:vdWq})}

Let $\KK$ be a field and $f(x)\in \KK[x]$ be a polynomial, not
necessarily irreducible, over $\KK$ of degree~$n$
with leading coefficient $=1${}.
A splitting field $\NN = \text{SF}_{\KK}(f)$
of $f$ over $\KK$ is constructed by a tower of
fields
\begin{equation}\label{eq:tower}
   \LL_0 = \KK \subset \LL_1 \subset \LL_2 \subset \cdots \subset \LL_\ell = \NN,
   \quad \LL_i = \LL_{i-1}[y_i]/(g_i(y_i)),\quad i=1,2,\ldots, \ell,
\end{equation}
where $y_1,\ldots,y_\ell$ are fresh variables and
where $g_i(x) \in \LL_{i-1}[x]$ is
an irreducible factor in $\LL_{i-1}[x]$ with $\deg_x(g_i) \ge 2$ of
\begin{equation}
f(x)\text{ if $i = 1$\quad and\quad}
\frac{f(x)}{(x-y_1)\ldots(x-y_{i-1})} \in \LL_{i-1}[x] \text{ if $i\ge 2$}.
\end{equation}
At index $\ell$ we have $f(x) = (x-\alpha_1)\cdots(x-\alpha_n)$ for $\alpha_i\in\NN${}.
All fields $\NN$ constructed in the manner (\ref{eq:tower}) are isomporhic with an
isomorphism that is the identity function on $\KK${}.
Note that for all $1 \le i \le \ell$
the fields $\LL_i$ are the quotient rings $\KK[y_1,\ldots,y_i]$ modulo the
triangular set $g_1(y_1),g_2(y_1,y_2),\ldots,g_i(y_1,\ldots,y_i)$
over $\KK${}.
Arithmetic in $\LL_i$ is done recursively as univariate polynomial residue
arithmetic in $\LL_{i-1}[y_i]/(g_i(y_i))$.
By (\ref{eq:tower}), $\NN = \text{SF}_{\LL_i}(f)$ for all $0 \le i \le \ell${}.
The Galois group
$\Gamma_{\NN/\KK}$ of $f(x)$ over $\KK$ is the group
of all field automorphisms $\psi\colon \NN \longrightarrow \NN$ with
$\psi(a) = a$ for all $a \in \KK${}.
Each automorphism $\psi$ uniquely permutes the distinct roots of $f$:
$\psi(\alpha_i) = \alpha_{\tau(i)}$,
and if $f(x)$ is separable, which means all roots are distinct:
$\alpha_i \ne \alpha_j$ for all $1 \le i < j \le n$, then
$\tau \in S_n$, where $S_n$ is the symmetric group of permutations on $1,\ldots,n$,
and all permutations $\tau$ form a subgroup.

%omittext

In \citeP[Section~61]{vdWae40} the permutations $\tau$
in the Galois group of a separable polynomial $f$ over a field $\KK$
are characterized as follows.
\begin{theorem}\label{thm:galois}
Let $f(x) = \prod_{i=1}^n (x-\alpha_i) \in \KK[x]$
where $\alpha_i\in\NN = \text{\upshape SF}_\KK(f)$
with $\alpha_i \ne \alpha_j$ for all $1 \le i < j \le n$,
and let
\begin{equation}\label{eq:Fzu}
F(z,u_1,\ldots,u_n) = \prod_{\sigma\in S_n}
\Big(z - \big(\sum_{i=1}^n \alpha_{\sigma(i)} u_i\big)\Big)
\in \KK[z,u_1,\ldots,u_n].
\end{equation}
Furthermore, let $F_1$ be an irreducible factor of $F$ in\/ $\KK[z,u_1,\ldots,u_n]$
such that $z - (\sum_{i=1}^n \alpha_i u_i)$ is a factor of $F_1$ over $\NN${}.
Then the permutations $\tau$ in the Galois group of $f$ over\/ $\KK$
are exactly those permutations such that
$z - (\sum_{i=1}^n \alpha_{\tau(i)} u_i)$ is a factor of $F_1$.
\end{theorem}
\ifomittext\else
\begin{proof}
Let $F_1 = \big(z - (\sum_{i=1}^n \alpha_i u_i)\big) G_1$ with
$G_1\in\NN[z,u_1,\ldots,u_n]${}.
Then $F_1 = \psi(F_1) = \big(z - (\sum_{i=1}^n \psi(\alpha_i) u_i)\big)\;
\psi(G_1)$ for all $\psi\in\Gamma_{\NN/\KK}${}.
Because $f$ is separable
all $\sum_{i=1}^n \psi(\alpha_i) u_i$ are distinct, and therefore all
$z - (\sum_{i=1}^n \psi(\alpha_i) u_i)$ divide $F_1$ over\/ $\NN$,
whose product is the norm in the splitting field
$\NN(z,u_1,\ldots,u_n)$ of $f(x)$ over $\KK(z,u_1,\ldots,u_n)$,
and therefore $\in \KK[z,u_1,\ldots,u_n]$.
\end{proof}
\fi%omittext
%In fact, $F_1$ is the norm of $z - (\sum_{i=1}^n \alpha_i u_i)$
%in the splitting field of $f$ over $\KK$ and as such a power of an irreducible
%polynomial in $\KK[z,u_1,\ldots,u_n]$, but which is also a separable polynomial
%in~$z$.

Note that the assumption that the roots $\alpha_i$ of $f$ are distinct
is a necessary condition.  Let $\KK = \FF 2(t)$ and $f(x) = x^2+t = (x+\sqrt{t})^2$.
Then $F(z) = (z + u_1 \sqrt{t} + u_2 \sqrt{t})^2 = z^2 + u_1^2 t + u_2^2 t$,
which is irreducible over $\FF 2(t)[z,u_1,u_2]$, but the Galois group
of $f(x)$ over $\FF 2(t)$ has a single element.
Because $x^2+t$ is irreducible in $\FF q(t)[x]$, it is squarefree
in $\FF q[x,t]$ but not squarefree (inseparable) over the algebraic closure
of $\FF q(t)$.

For generic polynomials the Galois group is the full symmetric group for all fields.

\begin{theorem}\label{thm:generic}
For the generic polynomial $f^{[v]} = x^n + \sum_{i=0}^{n-1} v_i\, x^i$
over\/ $\KK^{[v]} = \KK(v_0,\ldots,v_{n-1})$ the
polynomial $F^{[v]}$ corresponding to {\upshape (\ref{eq:Fzu})}
is a separable polynomial in $z$,
hence $\partial F^{[v]}/\partial z \ne 0$,
and an irreducible polynomial in
$\KK[z,u_1,\ldots,u_n,v_0,\ldots,v_{n-1}]$, for all fields~$\KK$.
\end{theorem}
\ifomittext\else
\begin{proof}
First, $f^{[v]}(x)$ is separable in~$x$
because it is irreducible over $\KK^{[v]}$
and its derivative with respect to $z$ is $\ne 0$.
The univariate polynomial discriminant is a non-zero polynomial
in the coefficients over
fields of all characteristics, which is $\ne 0$ for exactly the separable
polynomials.  Therefore, $F^{[v]}$ is also separable in $z$
implying that $\partial F^{[y]}/\partial z \ne 0$.

Let $\prod_{i=1}^n (x-w_i)
= x^n + e_{n-1}(w_1,\ldots,w_n)x^{n-1}+\cdots+e_0(w_1,\ldots,w_n)
\in \KK[z$, $w_1,\ldots$, $w_n]$,
where $e_i$ are plus/minus the $(n-i)$'th elementary
symmetric functions in fresh variables $w_1,\ldots,w_n$,
and let $\bar F^{[v]}$ be $F^{[v]}$ evaluated at $v_i =  e_i(w_1,\ldots,w_n)${}.
We have
$\bar F^{[v]} = \prod_{\sigma\in S_n}\big(z - \sum_{i=1}^n w_{\sigma(i)} u_i)${}.
Now let $F_1^{[v]}$ be an irreducible factor of $F^{[v]}$
in $\KK[z,u_1,\ldots,u_n$, $v_0,\ldots$, $v_{n-1}]$
and let $\bar F_1$ be $F_1^{[v]}$
evaluated at $v_i = e_i(w_1,\ldots,w_n)${}.
Then by definition of $\bar F^{[v]}$, there is a permutation $\tau\in S_n$ such
that $z - (w_{\tau(1)} u_1 + \cdots + w_{\tau(n)} u_n)$ divides $\bar F_1$ with
co-factor $\bar G_1\in\KK[z,u_1,\ldots,u_n,w_1,\ldots,w_n]${}.
Permuting the $w_i$'s in that factorization of $\bar F_1$  does not
change $\bar F_1$ and shows that
$z - (w_{\sigma(1)} u_1 + \cdots + w_{\sigma(n)} u_n)$ divides $\bar F_1$
for all permutations $\sigma\in S_n$.  Therefore $F_1$ has degree $n!$ in $z$.
\end{proof}
\fi%omittext

%%% roots of $f(x) = \prod_{i=1}^n (x - \alpha_i^{[y]})$, hence
%%% with coefficients in $\KK[y_0,\ldots,y_{n-1}]$.
%%% Now let $F_1^{[y]}$ be an irreducible factor of $F^{[y]}$
%%% over $\KK^{[y]}$ that has
%%% $z - (\sum_{i=1}^n \alpha_i^{[y]} u_i)$ as a factor over the
%%% algebraic closure of $\KK^{[y]}$.
%%% We can write
%%% $y_i = (-1)^{n-i} e_{n-i}(\alpha_1^{[y]},\ldots,\alpha_n^{[y]})
%%% \in \KK[\alpha_1^{[y]},\ldots,\alpha_n^{[y]}]
%%% $,
%%% where
%%% $e_{n-i}$ are the elementary symmetric functions.
%%% Then
%%% \begin{equation}\label{eq:facFone}
%%% F_1^{[y]}=\big(z - (\sum_{i=1}^n \alpha_i^{[y]} u_i)\big)
%%% \; G_1(z,u_1,\ldots,u_n) \in
%%% \KK[\alpha_1^{[y]},\ldots,\alpha_n^{[y]}][z,u_1,\ldots,u_n].
%%% \end{equation}
%%% Permuting the $\alpha_i^{[y]}$ we obtain the same left side in (\ref{eq:facFone}),
%%% %It follows that any irreducible factor $F_1^{[y]}$ over $\KK^{[y]}$ that has
%%% %$z - (\sum_{i=1}^n \alpha_i^{[y]} u_i)$ as a factor over the
%%% %algebraic closure of $\KK^{[y]}$
%%% implying that all
%%% $z - (\sum_{i=1}^n \alpha_i^{[y]} u_{\sigma(i)})$ are factors
%%% of $F_1^{[y]}$ over the algebraic closure of $\KK^{[y]}$,
%%% and therefore $F_1^{[y]} = F^{[y]}$, which concludes the
%%% proof of irreducibility of $F^{[y]}$ over $\KK^{[y]}$.
%%% %In conclusion $F^{[y]}$, whose leading coefficient in $z$ is $=1$,
%%% %is an irreducible polynomial in $\KK[z,u_1,\ldots,u_n,y_0,\ldots,y_{n-1}]${}.

Classically, one uses the Hilbert Irreducibility Theorem to count
for which evaluations of the $v_i$ at values in $\KK$ one preserves
irreducibility of $F$ \citeP{Knob56}.
For $\KK = \FF q(t)$ we can use our effective Hilbert Irreducibility
Theorems \citeP{Ka85:infcontr,Ka95:jcss}.  We have the following theorem.

\begin{theorem}
\label{thm:hit}
Let $F(X_1,\ldots,X_m)\in\KK[X_1,\ldots,X_m]$,
$\KK$ a field, have total degree $\delta$ and be irreducible.
Assume that $\partial F/\partial X_m \ne 0$.
%If the characteristic of $\KK$ is $p \ge 2$,
%we require that each coefficient of $F$
%in $\KK$ possesses a $p$-th root in $\KK$.
%A sufficient condition for this to be true is that $\KK$ is perfect.
Let $S\subseteq \KK$ be a finite set,
and let $a_2,\ldots,$ $a_{m-1}$,
$b_1,\ldots,$ $b_{m-1}$ be randomly and uniformly sampled elements in $S$.
Then the probability
\begin{align}
&\text{\upshape Prob}\Big(
F(b_1 , b_2,\ldots, b_{m-1} , z ) \in \KK[z]
\text{ is of degree }\deg_{X_m}(F)
\text{ and has discriminant $\ne 0$}
\notag
\\
&\text{and }
F(t +b_1 , a_2 t + b_2,\ldots,a_{m-1} t + b_{m-1} , z )
\text{ is irreducible in }\KK[t , z]\Big)
\ge 1 - \frac{4 \delta~2^\delta}{|S|},
\label{eq:hit}
\end{align}
where $|S|$ is the number of elements in the set $S$
{\upshape \citeP[Theorem 2 and its proof]{Ka85:infcontr}}.
\end{theorem}

We apply Theorem~\ref{thm:hit} to 
\begin{equation}
F^{[v]}(z, u_1,\ldots,u_n, v_0,\ldots,v_{n-1})
\in\KK(u_1,\ldots,u_n)[z, v_0,\ldots,v_{n-1}],
\end{equation}
which is defined above for the generic $f^{[v]}(x)$.
The leading coefficient of $F^{[v]}$ in $z$ is $=1$
and $F^{[v]}$ is irreducible over $\KK(u_1,\ldots,u_n)${}.
%%% Note that %$\partial F^{[y]}/\partial z \ne 0$
%%% the discriminant of $F^{[y]}$ in the variable $z$ is $\ne 0$
%%% because all $\sum_{i=1}^n  \alpha_i^{[y]} u_{\sigma(i)}$ are distinct.
%%% Therefore $\partial F^{[y]}/\partial z \ne 0$.
We have for randomly and uniformly sampled
$a_1,\ldots,a_{n-1}$,
$b_0,\ldots,b_{n-1} \in S \subseteq \KK \subset \KK(u_1,\ldots,u_n)$
and
\begin{equation}
\overline{F^{[v]}}(z,u_1,\ldots,u_n,t)
\defequal F^{[v]}(z,u_1,\ldots,u_n,t + b_0,a_1 t + b_1,\ldots,a_{n-1} t +b_{n-1})
\end{equation}
the probability estimate
\begin{align}
\text{\upshape Prob}\Big( &
\text{the discriminant of }
\overline{F^{[v]}}(z,u_1,\ldots,u_n,t)
\text{ in the variable $z$ is $\ne 0$ and}
\notag
\\
&\overline{F^{[v]}}(z,u_1,\ldots,u_n,t)
\text{ is irreducible in }\KK[z,u_1,\ldots,u_n,t]\Big)
\ge 1 - \frac{4 \delta^{[v]}~2^{\delta^{[v]}}}{|S|},
\label{eq:vdWbar}
\end{align}
where $\delta^{[v]}$ is the total degree of ${F^{[v]}}$
in $z,v_0,\ldots,v_{n-1}$.
% We argue by contradiction.  Suppose the fraction
% % $N/|S|^{2n}$
% of polynomials
% $\overline{F^{[y]}}(z,u_1,\ldots,u_n,t)$ for which the discriminant
% in $z$ is $=0$ or which are reducible in $\KK[z,u_1,\ldots,u_n,t]$
% is $> {4 \delta^{[y]}~2^{\delta^{[y]}}}/{|S|}$.
% Then for all evaluations of $u_1 = t + b_n$,
% $u_2 = a_n t + b_{n+1},\ldots, u_n = a_{2n-2} t + b_{2n-1}$
% with $a_i \in S$ for $n \le i \le 2n-2$ and $b_i \in S$ for
% $n \le i \le 2n-1$ the corresponding evaluated polynomials
% $\overline{F^{[y]}}(z,t + b_n,a_n t + b_{n+1},\ldots, a_{2n-2} t + b_{2n-1},t)
% \in \KK[z,t]$ are reducible or their discriminants with respect to the
% variable $z$ are $=0${}.
% By Theorem~\ref{thm:hit} those are too many
% polynomials.
All polynomials
$\bar f(x) = x^n + (a_{n-1} t + b_{n-1}) x^{n-1} + \cdots + (a_0 t + b_0)$
for which
$\overline{F^{[v]}}(z,u_1,\ldots,u_n,t)$ is irreducible and separable,
the latter of which implies that $\bar f(x)$ is separable,
have Galois group $S_n$ over $\KK(t)${}.
Because $n$ is a constant, $\delta^{[v]}$ is a constant.  For
the actual probability estimate (\ref{eq:vdWq}) we can set $\KK = S = \FF q$
and $t = t/a_0$ and multiply the probability~(\ref{eq:vdWbar})
by $(1 - 1/|S|)$ for $a_0 \ne 0${}.
The more specific evaluation $X_1 = t + b_1$ in Theorem~\ref{thm:hit}
strengthens our effective Hilbert Irreducibility Theorem.

\section{Remarks}
\label{sec:rems}

Better estimates than (\ref{eq:hit}) in terms of the degree for the effective Hilbert
Irreduciblity Theorems for function fields are possible.
An estimate $1 - O(\deg(F)^4/|S|)$ is in \citeP{Ka95:jcss} for perfect fields $\KK$,
which includes all $\FF q$.

The estimate (\ref{eq:vdWt}) follows from (\ref{eq:vdW})
by counting the irreducible $F^{[y]}(z,u_1,\ldots,u_n)$
for $y_i = a_i + t b_i$ with integers bounded by $|a_i| \le H$ and $|b_i| \le H$ and
the variable evaluation $t = 2H+1$ which implies $|y_i| \le 2H^2 + 2H$,
with $(2H+1)^2$ values for each $y_i${}.
The count implies that
$\text{GCD}(x^n + \sum_{i=0}^{n-1} a_i x^i, \sum_{i=0}^{n-1} b_i x^i) \ne 1$
occurs with probability $O(1/H^2)$ for fixed $n$.

\vspace*{\bigskipamount}\noindent
{\bfseries Acknowledgment:}
I thank Theresa C. Anderson for her correspondence about the topic of the paper.

\def\refname{\Large\bfseries References}

\bibliographystyle{plainnat}

\begin{thebibliography}{7}
\expandafter\ifx\csname natexlab\endcsname\relax\def\natexlab#1{#1}\fi
\expandafter\ifx\csname url\endcsname\relax
  \def\url#1{{\tt #1}}\fi

\bibitem[Anderson et~al.(2021)Anderson, Gafni, Oliver, Lowry-Duda, Shakan, and
  Zhang]{AGOL-DSZ21}
Anderson, Theresa~C., Gafni, Ayla, Oliver, Robert J.~Lemke, Lowry-Duda, David,
  Shakan, George, and Zhang, Ruixiang.
\newblock Quantitative {Hilbert} irreducibility and almost prime values of
  polynomial discriminants, 2021.
\newblock URL: \url{https://arxiv.org/abs/2107.02914}.

\bibitem[Benjamin and Bennett(2007)]{BB07}
Benjamin, Arthur~T. and Bennett, Curtis~D.
\newblock The probability of relatively prime polynomials.
\newblock {\em Mathematics Magazine}, 80\penalty0 (3):\penalty0 196--202, 2007.
\newblock URL: \url{https://doi.org/10.1080/0025570X.2007.11953481}.

\bibitem[Bhargava(2021)]{bhargava2021galois}
Bhargava, Manjul.
\newblock Galois groups of random integer polynomials and van der {Waerden}'s
  {Conjecture}, 2021.
\newblock URL: \url{https://arxiv.org/abs/2111.06507}.

\bibitem[Kaltofen(1985)]{Ka85:infcontr}
Kaltofen, E.
\newblock Effective {Hilbert} irreducibility.
\newblock {\em Information and Control}, 66:\penalty0 123--137, 1985.
\newblock
  \EKhref{http://users.cs.duke.edu/~elk27/bibliography/85/Ka85_infcontr.pdf}
  {EKbib/85/Ka85_infcontr.pdf}.

\bibitem[Kaltofen(1995)]{Ka95:jcss}
Kaltofen, E.
\newblock Effective {Noether} irreducibility forms and applications.
\newblock {\em J. Comput. System Sci.}, 50\penalty0 (2):\penalty0 274--295,
  1995.
\newblock
  \EKhref{http://users.cs.duke.edu/~elk27/bibliography/95/Ka95_jcss.pdf}
  {EKbib/95/Ka95_jcss.pdf}.

\bibitem[Kobloch(1956)]{Knob56}
Kobloch, Hans-Wilhelm.
\newblock Die {Seltenheit} der reduziblen {Polynome}.
\newblock {\em Jahresbericht d. DMV}, 59:\penalty0 12--20, 1956.
\newblock URL:
  \url{http://resolver.sub.uni-goettingen.de/purl?PPN37721857X_0059}.

\bibitem[{van der Waerden}(1940)]{vdWae40}
{van der Waerden}, B.~L.
\newblock {\em {Moderne Algebra}}.
\newblock Springer Verlag, Berlin, 1940.
\newblock English transl. publ. under the title ``Modern algebra'' by F. Ungar
  Publ. Co., New York, 1953.

\end{thebibliography}

\newcommand{\Erich}{Erich L.\xspace} \renewcommand{\Erich}{Erich L. }
  \newcommand{\Gathen}{\relax} \newcommand{\Hoeij}{\relax} \newcommand{\ZZZ}{}

\section{Appendix} 
The norm
$\text{norm}_{\NN/\KK}(\beta(y_1,\ldots,y_\ell))$
of an element $\beta(y_1,\ldots,y_\ell)\in \NN$ over $\KK$,
where $\NN$ is the splitting field~(\ref{eq:tower})
of a possible inseparable polynomial $f$,
is defined recursively:
\begin{equation}\label{eq:norm}
\text{norm}_{\LL_{\ell-1}/\KK}\underbrace{\Big(
\prod_{j=1}^k
\beta(y_1,\ldots,y_{\ell-1},\gamma_j)\Big)}_{
\text{norm}_{\NN/\LL_{\ell-1}}(\beta(y_1,\ldots,y_\ell))\in\LL_{\ell-1}
\footnotemark
}
\in\KK
,
g_\ell(x) = (x-\gamma_1)\cdots(x-\gamma_k),
\gamma_i\in\NN, \gamma_1 = y_\ell.
\end{equation}
\footnotetext{%
Note that $\text{norm}_{\NN/\LL_{\ell-1}}(\beta(y_1,\ldots,y_\ell))$
is the Sylvester resultant of
$g_\ell(x)$ and $\beta(y_1,\ldots,y_{\ell-1},x)$
with respect to the variable $x$.
}%
The definition (\ref{eq:norm}) extends to the rational function fields
$\NN(X_1,\ldots,X_m)$ over $\KK(X_1,\ldots,X_m)$, where we have the following
theorem.
\begin{theorem}\label{thm:norm}
Let $G \in \NN[X_1,\ldots,X_m]$ be an irreducible polynomial
over\/ $\NN$, where $\NN$ is the splitting field\/ {\upshape(\ref{eq:tower})}
of a not necessarily separable polynomial.
Then $\text{norm}_{\NN/\KK}(G) = H^k$ where
$H\in\KK[X_1,\ldots,X_m]$ is irreducible over\/ $\KK$ and $k\ge 1$.
\end{theorem}
\begin{proof}
Suppose $\text{norm}_{\NN/\KK}(G) = H_1 H_2$ with
$H_1,H_2\in \KK[X_1,\ldots,X_m]$ and
$\text{GCD}(H_1,H_2) = 1$.  Note that relatively primeness
as an arithmetic property over $\KK$ remains valid over $\NN${}.
%Let $\bar H_1$ be an irreducible factor of $H_1$ in $\KK[X_1,\ldots,X_m]$, and
Now suppose
that $G(X_1,\ldots,X_m,y_1,\ldots,y_\ell)$ is an irreducible factor of $H_1$
over $\NN$.
By definition (\ref{eq:norm})
there exist roots $\gamma_i \in \NN$ of $g_i(x)$
such that $G(X_1,\ldots,X_m,\gamma_1,\ldots,\gamma_\ell)$
divides $H_2$ over $\NN${}.
The field $\NN$ is isomorphic to $\KK(\gamma_1,\ldots,\gamma_\ell)$
by $\psi\colon y_i \mapsto \gamma_i$ and $\psi(a) = a$ for all $a\in\KK$,
so $G(X_1,\ldots,X_m,\gamma_1,\ldots,\gamma_\ell)$ divides $\psi(H_1) = H_1$
over $\NN$, which contradicts that $H_1,H_2$ are relatively prime.
\end{proof}

Note that for $\beta\in\NN$
we have $\text{norm}_{\NN/\KK}(x-\beta) = h(x)^k$
where $h(x)\in\KK[x]$ is the irreducible minimum polynomial
with $h(\beta)=0$, which means that $\text{norm}_{\NN/\KK}(\beta)$
is $k$-th power of the product of all conjugates of $\beta$ over~$\KK$, which
are the roots of $h$ with multiplicities.
For a separable polynomial $f(x)$ and $\beta\in\NN=\text{SF}_\KK(f)$,
we have $\text{norm}_{\NN/\KK}(\beta)
= \prod_{\psi\in\Gamma_{\NN/\KK}} \psi(\beta)$, where
$\Gamma_{\NN/\KK}$ is the Galois group as a group of field automorphisms.

\vspace{\bigskipamount}\noindent{\itshape Proof of Theorem \ref{thm:galois}.}
Let $F_1 = \big(z - (\sum_{i=1}^n \alpha_i u_i)\big) G_1$ with
$G_1\in\NN[z,u_1,\ldots,u_n]${}.
Then $F_1 = \psi(F_1) = \big(z - (\sum_{i=1}^n \psi(\alpha_i) u_i)\big)\;
\psi(G_1)$ for all $\psi\in\Gamma_{\NN/\KK}${}.
Because $f$ is separable
all $\sum_{i=1}^n \psi(\alpha_i) u_i$ are distinct, and therefore all
$z - (\sum_{i=1}^n \psi(\alpha_i) u_i)$ divide $F_1$ over\/ $\NN$,
whose product is the norm in the splitting field
$\NN(z,u_1,\ldots,u_n)$ of $f(x)$ over $\KK(z,u_1,\ldots,u_n)$,
and therefore $\in \KK[z,u_1,\ldots,u_n]$.~$\Box$

\vspace{\bigskipamount}\noindent{\itshape Second proof of Theorem \ref{thm:galois}.}
By Theorem~\ref{thm:norm}
the norm of $z - \sum_{i=1}^n \alpha_i u_i$ is $H(z,u_1,\ldots,u_n)^k$ with
$H$ irreducible in $\KK[(z,u_1,\ldots,u_n]${}.  The norm's discriminant in $z$
is $\ne 0$ because the roots are distinct, which implies $k=1$.~$\Box$

\vspace{\bigskipamount}\noindent{\itshape Proof of Theorem \ref{thm:generic}.}
First, $f^{[v]}(x)$ is separable in~$x$
because it is irreducible over $\KK^{[v]}$
and its derivative with respect to $z$ is $\ne 0$.
The univariate polynomial discriminant is a non-zero polynomial
in the coefficients over
fields of all characteristics, which is $\ne 0$ for exactly the separable
polynomials.  Therefore, $F^{[v]}$ is also separable in $z$
implying that $\partial F^{[v]}/\partial z \ne 0$.

Let $\prod_{i=1}^n (x-w_i)
= x^n + e_{n-1}(w_1,\ldots,w_n)x^{n-1}+\cdots+e_0(w_1,\ldots,w_n)
\in \KK[z$, $w_1,\ldots$, $w_n]$,
where $e_i$ are plus/minus the $(n-i)$'th elementary
symmetric functions in fresh variables $w_1,\ldots,w_n$,
and let $\bar F^{[v]}$ be $F^{[v]}$ evaluated at $v_i =  e_i(w_1,\ldots,w_n)${}.
We have
$\bar F^{[v]} = \prod_{\sigma\in S_n}\big(z - \sum_{i=1}^n w_{\sigma(i)} u_i)${}.
Now let $F_1^{[v]}$ be an irreducible factor of $F^{[v]}$
in $\KK[z,u_1,\ldots,u_n$, $v_0,\ldots$, $v_{n-1}]$
and let $\bar F_1$ be $F_1^{[v]}$
evaluated at $v_i = e_i(w_1,\ldots,w_n)${}.
Then by definition of $\bar F^{[v]}$, there is a permutation $\tau\in S_n$ such
that $z - (w_{\tau(1)} u_1 + \cdots + w_{\tau(n)} u_n)$ divides $\bar F_1$ with
co-factor $\bar G_1\in\KK[z,u_1,\ldots,u_n,w_1,\ldots,w_n]${}.
Permuting the $w_i$'s in that factorization of $\bar F_1$  does not
change $\bar F_1$ and shows that
$z - (w_{\sigma(1)} u_1 + \cdots + w_{\sigma(n)} u_n)$ divides $\bar F_1$
for all permutations $\sigma\in S_n$.  Therefore $F_1$ has degree $n!$ in $z$.~$\Box$

\end{document}